\documentclass[12pt]{article}
\usepackage{amsmath,amssymb,amsthm}

\numberwithin{equation}{section}

\newtheorem{thm}{Theorem}[section]
\newtheorem{lem}{Lemma}[section]

\newcommand{\n}{\nonumber}

\newcommand{\si}{\sigma_R }

\renewcommand{\a}{\alpha}
\renewcommand{\l}{\lambda}
\newcommand{\vare}{\varepsilon}
\newcommand{\s}{\sigma}
\renewcommand{\o}{\omega}

\renewcommand{\O}{\Omega}
\newcommand{\vb}{\overline{V}}
\newcommand{\Ob}{\overline{\Omega}}
\newcommand{\se}{(SSE)_{\alpha}}
\newcommand{\bb}{\begin{equation}}
\newcommand{\ee}{\end{equation}}
\newcommand{\bq}{\begin{eqnarray}}
\newcommand{\eq}{\end{eqnarray}}
\newcommand{\bqn}{\begin{eqnarray*}}
\newcommand{\eqn}{\end{eqnarray*}}
\begin{document}
\title{ Euler's equations and the maximum principle}
\author{Dongho Chae \\
 Department of Mathematics\\
Chung-Ang University\\
 Seoul 156-756, Republic of Korea\\
e-mail: dchae@cau.ac.kr}
\date{(to appear in Math. Ann.)}
\maketitle

\begin{abstract}
In this paper we use maximum principle in the far field region for the time dependent self-similar Euler equations to exclude  discretely self-similar blow-up for the Euler equations of the incompressible fluid flows. Our decay conditions near spatial infinity
of the blow-up profile are given explicitly in terms the coefficient in the equations.  We also deduce triviality of the  discretely self-similar solution to the magnetohydrodynamic system(MHD), under suitable decay conditions near spatial infinity than the previous one. Applying similar argument directly to the Euler equations, we obtain a priori estimate of the vorticity in the far field region. \\
\ \\
\noindent{\bf AMS Subject Classification Number:} 35Q31, 76B03,
76W05\\
  \noindent{\bf
keywords:}   Euler equations, discretely self-similar solution, maximum principle, MHD
\end{abstract}
\section{The main theorems}
\setcounter{equation}{0}
We are concerned on the homogeneous incompressible 3D Euler equations.
$$
(E) \left\{\aligned  &\frac{\partial v}{\partial t }+v\cdot \nabla v =-\nabla p,\label{e1}\\
&\mathrm{div }\, v=0,\\
& v(x,0)=v_0(x),
\endaligned \right.
$$
 where $v(x,t)=(v_1 (x,t), v_2 (x,t), v_3 (x,t))$ is the velocity, $p=p(x,t)$ is the
  pressure, and $v_0(x)$ is the initial data satisfying div$v_0 =0$. For the Cauchy problem of (E) the  local well-posedness is known  for  the initial data belonging to the standard Sobolev space, $v_0 \in H^m (\Bbb R^3), m>5/2$(see e.g. \cite{kat}). The question of spontaneous apparition of singularity in finite time, however,  is still an outstanding open problem(see \cite{bea,con-fef}). For a general introduction to the blow-up problem for smooth initial data  we refer  a book(e.g. \cite{maj}) or surveys(\cite{con, bar}).
  Here we concentrate on the study of specific scenario of blow-up of self-similar type, and its generalized version, the discretely self-similar blow-up(see the definition below). They have been studied previously in a series of papers(e.g. \cite{cha1,cha2, cha3, cha-shv, cha-tsai}). An important question in this direction of study is how much we can relax the decay condition near spatial infinity of the blow-up profile and deduce its triviality, hence leading to  exclude  self-similar type blow-up in a full generality. The similar question for the Navier-Stokes equations was already answered satisfactorily at the level of the self-similar blow-up (not yet for the discretely self-similar blow-up) in the works of \cite{nec, tsai, mil}. See also \cite{egg, gig1, gig2} for the earlier studies of the self-similar solutions in other nonlinear partial differential equations. For a comprehensive introduction to the self-similar solution with its various aspects we refer the book \cite{gig0}.

  In this paper, developing new argument of using the maximum principle in the far field region, we deduce the triviality of the discretely self-similar solutions to the system (E), under  mild decay conditions near spatial infinity for the blow-up profile, thus substantially improving  and complementing  the previously known results. This type of argument is completely  new in the study of
  self-similar solutions of the Euler equations as far as the author knows. We also show that this argument could be used for (E) directly to derive a priori estimate of the vorticity in the far field region. In the last section we show that our method of proof can also be applied to the exclusion of discretely self-similar solutions of more general system of inviscid magnetohydrodynamics, which improves also the previous results in this direction(\cite{cha3}).  \\

  Below, after reviewing briefly the notion of the discretely self-similar solution, we state our main results of the paper.
Given $\alpha \neq -1$, we make the  self-similar
transform of (E), which is defined by the map $(v,p) \mapsto (V,P)$
given by
\begin{equation}
\label{11}
 v(x, t)=\frac {1}{(T-t)^{\frac{\a}{\a+1}}} V(y,s), \quad
 p(x,t)=\frac {1}{(T-t)^{\frac{2\a}{\a+1}}} P(y,s),
\end{equation}
where
\begin{equation}
\label{12} y = \frac{x-x_*}{(T-t)^{\frac{1}{\a+1}}}, \quad s =  \log \left(\frac{T}{T-t}\right).
\end{equation}
Substituting (\ref{11})--(\ref{12}) into (E), we have the following system in terms of  $(V,P)$:
\begin{equation}\label{} (SSE)_\a \left\{
\aligned
&\frac {\partial V}{\partial s}+ \frac{\a}{\a+1} V +\frac{1}{\a+1}(y \cdot \nabla)V + (V\cdot \nabla )V =-\nabla
P,\\
& \mathrm{div}\,V=0,\\
& V(y,0)=V_0 (y)=T^{\frac{\a}{\a+1}}v_0 (T^{\frac{1}{\a+1}}y).
\endaligned \right.
\end{equation}
The self-similar solution $(v,p)$ of (E) is defined as a solution of (E) given by (\ref{11})-(\ref{12}), where $(V,P)$ is a {\em stationary} solution of $(SSE)_\a$. In the case when $T$ is the blow-up time of the solution $(v,p)$, then we say that the solution given by (\ref{11})-(\ref{12}) is a self-similar blowing-up solution.
This is a solution of (E), having  the scale symmetry with respect to $(x_*, T)$. Namely, after translation of the coordinate origin into $(x_*, T)$ we have the invariance,
\bq\label{13}
\l ^\a v(\lambda x, \lambda ^{\a+1} (-t))&=&v(x, -t), \n \\
 \l^{2\a} p(\lambda x, \lambda ^{\a+1} (-t))&=&p(x, -t)
\eq
{\em for all} $\l\in \Bbb R\setminus\{0\}$ and for all $(x,t)\in \Bbb R^3 \times (-\infty, 0)$.
On the other hand, there exists a weaker notion of the self-similarity for the solution, called discrete self-similarity.
We say a solution $(v,p)$ of (E) is a discretely self-similar solution to (E)  (with respect to $(x_*,T)$) if {\em there exists} $\lambda \neq 1$ such that (\ref{13}) holds true after translation of the origin.  We find that $(v,p)$ given by  (\ref{11})-(\ref{12}) is a discrete self-similar solution to (E) with $\l \neq 1$ if and only if the {\em time-dependent}  $(V,P)$ of $\se$ satisfies the periodicity in time,
\bb\label{period}
V(y, s)=V(y, s+S_0), \quad P(y,s)=P(y, s+S_0)\quad \mbox{with}\quad S_0=-(\a+1)\log \l , \quad \forall y\in \Bbb R^3.
\ee
Conversely, for any time periodic solution  $(V,P)$ of $\se$ with the period $S_0\neq 0$ one can generate discretely self-similar solution of (E) with the scaling parameter given by (\ref{period}).
Thus, the question of the existence of nontrivial discretely self-similar solutions is equivalent to that of the existence of nontrivial time-periodic solution to $\se$.
The following is our main result in this question. For $k,m\in \Bbb N\cup\{0\}$ we denote that $f=f(y,s)\in C^k_sC_y^m (\Bbb R^{3+1})$ if the function $s\mapsto f(y,s) $ is of the $C^k (\Bbb R)$ class for each $y\in \Bbb R^3$, and the function $y\mapsto f(y,s) $ is of the $C^m (\Bbb R)$ class for each $s\in \Bbb R$.
\begin{thm}
Let $V\in C^1_sC^2_y (\Bbb R^{3+1})$ be a time periodic solution to $\se$ with the period $S_0 \neq 0$,
which satisfies
  \bb\label{decay}
 |\nabla V( y,s)| =o(1) \quad \mbox{as $|y|\to \infty$}\quad \forall s\in [0, S_0].
  \ee
  In the case $\a <-1$ we do not assume any extra condition, while
 if $\a>-1 $, we assume that there exists $k>\a+1$ such that the vorticity $\O=\mathrm{curl} \, V$ satisfies
\bb\label{decayiii}
\ |\O(y,s)|=O (|y|^{-k}) \quad \mbox{as $|y|\to \infty$}\quad \forall s\in [0, S_0].
\ee
Then,  $V(y,s)=C(s)$ for all $y\in \Bbb R^3$, where $C:[0, S_0] \to \Bbb R^3$ is a closed curve, $C(s)=C(s+S_0)$ for all $s\in [0, s_0)$.
\end{thm}
\noindent{\em Remark 1.1 } In Theorem 4.1 of \cite{cha-shv} it is proved that for $\a >-1$  the {\em stationary solution}
to $\se$ is trivial, namely  the velocity $V$ is  a constant vector, if $\O\in L^p(\Bbb R^3)$ with $p<\frac{3}{\a+1}$ together with the condition (\ref{decay}), which is a comparable condition with the above.
In the case $\a <-1 $, however,  the above result is completely new. Similar remarks hold, comparing the above theorem with  Theorem 2.2 of \cite{cha-tsai} for the time-periodic solutions of $\se$.\\
\ \\
Let us consider the time-periodic solutions of the following more general system than $\se$.
\bb\label{v}
 \left\{\aligned  &\frac{\partial V}{\partial s}  +aV+b (y\cdot\nabla ) V+ (V\cdot \nabla) V =-\nabla P, \\
&\mathrm{div }\, V=0,\\
&V(y,0)=V_0 (y),
\endaligned \right.
\ee
 where $a,b\in \Bbb R$, and $b\neq 0$. For the system (\ref{v}) we have the following result, from which  Theorem 1.1 follows as an immediate corollary.
\begin{thm}
Let $V\in C_s ^1C^2_y(\Bbb R^{3+1})$ be a time periodic solution to (\ref{v}) with the period $S_0 \neq 0$,
  satisfying (\ref{decay}). In the case  $(a+b)b <0$ we assume no extra condition, while if $(a+b)b >0$
we assume that there exists $k>\frac{a+b}{b}$ such that
\bb\label{decayiia}
 \quad |\O(y,s)|=O(|y|^{-k}) \quad \mbox{as $|y|\to \infty$}\quad  \forall s\in [0, S_0].
\ee
Then, $V(y,s)=C(s)$ for all $y\in \Bbb R^3$, where $C:[0, S_0] \to \Bbb R^3$ is a closed curve, $C(s)=C(s+S_0)$ for all $s\in [0, s_0)$.
\end{thm}
The proof of Theorem 1.2 uses crucially the maximum principle in the far field region, restricted by suitable cut-off function, where the nonlinear vortex stretching term is dominated by the linear terms  under the hypothesis (\ref{decay}). The decay rate of the condition (\ref{decayiii}) is necessary
to adjust the strength of the linear terms in order to beat the nonlinearity.
We note that the crucial decay assumption (\ref{decay}), which seemingly  strong, is in fact redundant, and natural for the Cauchy problem of (E) with the standard initial data $v_0\in H^m (\Bbb R^3), m>5/2$, as  will be seen in the proof of Theorem 1.3 below.
The use of the maximum principle for the {\em inviscid problem} is
 the first time in this paper as far as the author knows.  As another application  of the argument to the Cauchy problem of (E) directly we can derive the following estimate for the vorticity in the far field region, which is similar to the 2D Euler equations.
\begin{thm}
Let $m>5/2$ and  $v\in C([0, T); H^m (\Bbb R^3))$ be a solution to (E) corresponding to the initial data $v_0\in H^m (\Bbb R^3)$,  where $T$  is the maximal time of local well-posedness. Then,
 for all $\vare>0$ and  $t_0\in (0, T)$ there exists $R_0=R_0 (v_0, t_0, \vare )$ such that
\bb\label{est}
 |\o(X_t(a),t)|\leq  e^{\vare t} |\o_0(a)|\qquad \forall (a,t)\in \left(\Bbb R^3\setminus B(0, R_0)\right) \times [0, t_0],
 \ee
where $\{ X_t (\cdot)\}$ is the particle trajectory mapping generated by the velocity field $v(x,t)$.
\end{thm}

\section{ Proof of the Main Theorems }
\setcounter{equation}{0}

Theorem 1.2 is an immediate consequence of combining the following two lemmas.
\begin{lem}
Let $V\in C_s ^1C^2_y(\Bbb R^{3+1})$ be a time periodic solution  with the period $S_0$ of (\ref{v})
with $b\neq 0$, satisfying (\ref{decay}). We assume that the decay conditions of Theorem 1.2  hold true.
Then, the vorticity $\O$ has compact support in $\Bbb  R^3\times [0, S_0]$.
\end{lem}
 \begin{lem}
Let $V\in C_s ^1C^2_y(\Bbb R^{3+1})$ be a time periodic solution  with the period $S_0$ of (\ref{v})
with $a\in \Bbb R$ and $b\neq 0$. If  we assume (\ref{decay}), and
\bb\label{comp}
\O \in  \bigcup_{r>0}\bigcap_{0<q<r} L^q(\Bbb R^{3}\times [0, S_0]),
\ee
then
$ V(y,s)=C(s)$ for all  $(y,s)\in \Bbb R^3\times[0, S_0]$, where $C: [0, S_0]\to \Bbb R^3$ is a closed curve with $C(s)=C(s+S_0)$ for all $s\in [0, s_0)$.
\end{lem}
\noindent{\em Remark 2.1 } We note the fact that, in particular, the condition (\ref{comp}) is  satisfied if  $\O $ has compact support in $\Bbb R^3\times [0, S_0]$ . Thus, combining Lemma 2.1 and Lemma 2.2, we are lead to the conclusion of Theorem 1.2.\\
\ \\
\noindent{\em Remark 2.2 } The proof of the {\em stationary version} of Lemma 2.2, having weaker conclusion under weaker hypothesis, $\nabla V\in L^\infty $ for the case $a=\frac{\a}{\a+1}, b=\frac{1}{\a+1}$ is done in \cite{cha2}, while its time-dependent version is  done in \cite{cha-tsai}.  Below we present the proof, which is applicable to our purpose, and has refined argument.\\
\ \\
\noindent{\bf Proof of Lemma 2.1 }  For $r>0$ we denote $B(0, r)=\{y\in \Bbb R^3\, |\, |y|<r\}$. We will prove it by contradiction. We assume that the support of $\O $ is non-compact
in $\Bbb R^3\times [0, S_0]$, namely for any large $R>0$ there exists an open set in $(\Bbb R^3\setminus B(0, R))\times [0, S_0]$,
where $|\O|> 0$. Then, we will derive contradiction depending on the decay conditions of $\O$ as in (i)-(iii).\\

\noindent{\em \underline{(A) The case of (i) with $a+b>0, b<0$:} }
Taking curl on (\ref{v}), we obtain the voticity equations,
\bb\label{vol}
\left\{\aligned &\frac {\partial \O}{\partial s}+ (a+b)\O + b (y\cdot \nabla) \O+(V\cdot \nabla )\O =(\O\cdot \nabla )V,\\
& \mathrm{curl} \,V=\O, \quad \mathrm{div} \, V=0,\\
& \O (y,s)=\O_0 (y).\endaligned \right.
\ee
For each $\rho, \delta >0$  we  define a cut-off function $\psi=\psi_\rho(x)$  as follows.
\bb\label{cut}
\psi_\rho (x)=\left\{ \aligned &1,  &\mbox{if}\quad |x|>\rho+\delta\\
&\frac12\sin\left(\frac{(|x|-\rho)\pi}{\delta}-\frac{\pi}{2}\right)+\frac12, &\mbox{if}\quad \rho<|x|\leq \rho+\delta\\
                                &0,        &\mbox{if} \quad|x|\leq \rho.
                                \endaligned
                                \right.
                                \ee
  Let $\delta >0$. Multiplying  (\ref{vol}) by $\O \psi_\rho $, we obtain
\bq\label{21}
\lefteqn{\frac{\partial}{\partial s} (\psi_\rho|\O|^2 ) + 2(a+b)\psi_\rho |\O|^2 +b(y\cdot \nabla )( \psi_\rho|\O|^2 ) +(V\cdot \nabla )(\psi_\rho|\O|^2 )}\hspace{.0in}\n \\
 &&=2\hat{A}  \psi_\rho|\O|^2+ b |\O|^2 (y\cdot \nabla )\psi_\rho +|\O|^2 (V\cdot \nabla )\psi_\rho.
 \eq
where we defined $\hat{A}(y,s)$ by
\bb\label{22} \hat{A}(y,s)=\left\{ \aligned & \frac{(\O(y,s))\cdot \nabla )V(y,s)\cdot \O(y,s)}{|\O(y,s)|^2}, \quad \mbox{if}\quad \O(y,s)\neq 0,\\
 & 0, \quad \mbox{if}\quad \O(y,s)= 0.
 \endaligned \right.
\ee
We note that
from the condition (\ref{decay})
we have the sub-linear growth at spatial infinity  for the radial component of the velocity, $V_r=V\cdot \frac{y}{|y|}$
\bb\label{23}
|V_r (y)|= o(|y|), \quad \mbox{as}\quad |y|\to \infty.
\ee
Indeed, we have
\bb\label{calculus}
\frac{V_r(y,s)}{|y|}= \frac{ V(y,s) \cdot y}{|y|^2}=  \frac{y}{|y|^2} \cdot \left\{V(0,s)+ \int_0 ^1 \nabla V (ty,s)\cdot y dt\right\}=o(1)
\ee
as $|y|\to \infty$ for all $s\in [0, S_0]$, and (\ref{23}) is verified.
From (\ref{decay})  one can choose sufficiently large $r_0=r_0(V_0)$  so that
 \bb\label{24}
|\hat{A}(y,s)|\psi_\rho(y) \leq |\nabla V(y,s)| \psi_\rho(y)
\leq \frac12(a+b) \psi_\rho (y),
 \ee
\bb\label{25}
 |(V\cdot \nabla )\psi_\rho|=  |V_r|\partial_r \psi_\rho \leq -\frac{b}{2} |y| \partial_r \psi_\rho
 =  -\frac{b}{2}   (y\cdot \nabla) \psi_\rho
\ee
 for all $\rho\geq r_0$ and for all $s\in [0, S_0]$.
Substituting (\ref{24}) and (\ref{25})  into (\ref{21}),  we  have
 \bq\label{26}
\lefteqn{\frac{\partial}{\partial s} (|\O|^2\psi_\rho) +2(a+b) |\O|^2\psi_\rho +b(y\cdot \nabla )(|\O|^2 \psi_\rho) +(V\cdot \nabla )(|\O|^2\psi_\rho)}\hspace{.1in}\n \\
&&=2\hat{A}(y,s) |\O|^2\psi_\rho(y)+b  |\O|^2(y\cdot \nabla )\psi_\rho +|\O|^2(V\cdot \nabla )\psi_\rho \n \\
&& \leq (a+b)  |\O|^2\psi_\rho (y)+ \frac{b}{2}  |\O|^2(y\cdot \nabla )\psi_\rho
 \eq
 for all $\rho \geq r_0$ and $s\in [0, S_0]$. Since $\psi_\rho$ is radially  non-decreasing, and $b<0$, we have $\frac{ b}{2} |\O|^2(y\cdot \nabla )\psi_\rho \leq 0$. Hence, we obtain the following differential inequality.
 \bb\label{27}
 \frac{\partial}{\partial s} f(y,s)+(a+b) f(y,s)+b y\cdot \nabla f(y,s)+V\cdot \nabla f(y,s)\leq 0,
 \ee
 where we set
 $$
 f(y,s):= |\O|^2\psi_{r_0}.
 $$
 Let us define the space-time domain
 $$\mathcal{D}_{r_0}:=(\Bbb R^3\setminus  B(0, r_0)) \times (0, S_0).$$
 We have $\partial \mathcal{D}_{r_0}= \Gamma_1 \bigcup \Gamma_2 \bigcup \Gamma_3 \bigcup \Gamma_4$, where
\bqn &&\Gamma_1=\partial B(0, r_0)\times  (0, S_0), \quad \Gamma_2=\{ |y|=\infty\}\times  (0, S_0) \\
&& \Gamma_3= (\Bbb R^3\setminus  B(0, r_0))   \times \{s= S_0\} ,\quad \Gamma_4= (\Bbb R^3\setminus B(0,r_0))\times \{ s=0\}.
\eqn
Since we assumed non-compactness of the support of $\O$ on $\Bbb R^3\times [0, S_0]$, there exists an open set in $\mathcal{D}_{r_0}$, where $f(y,s)$ is positive, which implies that there exists a positive maximum
of $f$ in $\mathcal{\overline{D}}_{r_0}$.
The differential inequality (\ref{27}) implies that the function $f(y,s)$ cannot have positive maximum in $\mathcal{D}_{r_0}\bigcup \Gamma_3$. Indeed, if $(\bar{y},\bar{s})\in \mathcal{D}_{r_0}\bigcup \Gamma_3$ is a point of the positive maximum for $f(y,s)$ on $\mathcal{D}_{r_0}\bigcup \Gamma_3$, then
$$
\lim_{s\uparrow \bar{s}} \frac{\partial}{\partial s}f(\bar{y},s)\geq 0,\quad  (y\cdot \nabla )f(\bar{y},\bar{s})=(V\cdot \nabla )f(\bar{y},\bar{s})=0,\quad
f(\bar{y},\bar{s})>0,
$$
which is a contradiction to (\ref{27}).
Since $\psi_{r_0}=0$ on $\Gamma_1$,  and $|\O (y,s)|\leq 2|\nabla V(y,s)|=0$ on $\Gamma_2$, we have
 $f(y,s)=0$ on $\Gamma_1 \bigcup \Gamma_2$.
Thus, the positive maximum of $f(y,s)$ on $\mathcal{\overline{D}}_{r_0}$ is attained only at $\Gamma_4$, and we have
 \bb\label{28}
 \sup_{|y|>r_0, s\in (0, s_0)}|\O(y,s)|^2\psi_{r_0}(y)\leq  \sup_{|y|>r_0}|\O_0(y)|^2\psi_{r_0}(y) .
\ee
Moreover,  there exists $y_0 \in \Bbb R^3$ with $|y_0|>r_0$ such that the strict inequality
\bb\label{29}
|\O(y,s)|^2\psi_{r_0}(y)< |\O_0(y_0)|^2\psi_{r_0}(y_0)=\sup_{|y|>r_0} |\O_0(y)|^2\psi_{r_0}(y)
\ee
 holds for all $(y,s)\in \mathcal{D}_{r_0}\bigcup \Gamma_3$. Substituting $y=y_0, s=S_0$ in (\ref{29}), we obtain
\bb\label{210}
|\O(y_0,S_0)|^2\psi_{r_0}(y_0)<  |\O_0(y_0)|^2\psi_{r_0}(y_0),
\ee
Since $\O(y_0,S_0)|=|\O_0(y_0)|>0$ by the periodicity, (\ref{210}) is absurd. \\
 \ \\

\noindent{\em \underline{(B) The case of (i) with $a+b<0, b>0$:} }
 In this case we define $\overline{V} (y, s)= V(y, S_0-s)$,  $\bar{P} (y, s)= P(y, S_0-s)$ and $\Ob (y, s)= \O(y, S_0-s)$ for $0\leq s\leq S_0$. Then, the vorticity equation becomes
 \bb\label{211}
\frac {\partial \Ob}{\partial s}-(a+b)\Ob - b (y\cdot \nabla) \Ob-(\vb\cdot \nabla )\Ob =-(\O\cdot \nabla )\vb.
\ee
 This is the same situation as (A) above with $a+b>0, b<0$. In particular, we note that the signs in front of the terms
 $(\overline{V}\cdot \nabla )\overline{\O}$ and $ (\overline{\O}\cdot \nabla )\overline{V}$ are not important in the estimates (\ref{24}), (\ref{25}).
 Hence, following the same argument as the proof (A), we conclude that $\Ob(y,s)$, and therefore $\O(y,s)$ has  compact support in  $\Bbb R^3\times [0, S_0]$.\\

 \noindent{\em \underline{(C) The case of (ii) with $a+b\leq 0, b<0 $:} } Let $k>\frac{a+b}{b}$.
 We  multiply  (\ref{vol}) by $\O \psi_\rho |y|^{2k} $ to obtain
\bq\label{212}
\lefteqn{\frac{\partial}{\partial s} (\psi_\rho|y|^{2k}|\O|^2 ) + 2(a+b-bk)\psi_\rho|y|^{2k} |\O|^2 +b(y\cdot \nabla )( \psi_\rho|y|^{2k}|\O|^2 ) +(V\cdot \nabla )(\psi_\rho|y|^{2k}|\O|^2 )}\hspace{.0in}\n \\
 &&\qquad=2\hat{A}  |y|^{2k}\psi_\rho|\O|^2+ b |y|^{2k}|\O|^2 (y\cdot \nabla )\psi_\rho +|y|^{2k}|\O|^2 (V\cdot \nabla )\psi_\rho +2k |y|^{2k-1} |\O|^2 V_r \psi_\rho.\n \\
 \eq
As previously
from the condition (\ref{decay})  one can choose $r_0=r_0(V_0)$  so that
 \bb\label{213}
|\hat{A}(y,s)|\psi_\rho(y) \leq \frac14 (a+b-bk)\psi_\rho (y),
 \ee
\bb\label{214}
 |(V\cdot \nabla )\psi_\rho|= |V_r|\partial_r \psi_\rho \leq -\frac{b}{2}  |y| \partial_r \psi_\rho
 =  -\frac{b}{2} (y\cdot \nabla) \psi_\rho,
\ee
and
\bb\label{214a}
2k |y|^{2k-1}  V_r \psi_\rho\leq 2k |y|^{2k} \psi_\rho \sup_{|y|>\rho} \frac{|V(y,s)|}{|y|}\leq \frac12 (a+b-bk) |y|^{2k} \psi_\rho
\ee
 for all $\rho\geq r_0$ and for all $s\in [0, S_0]$.
Substituting (\ref{213})-(\ref{214a})  into (\ref{212}),  we  have
 \bq\label{215}
\lefteqn{\frac{\partial}{\partial s} (\psi_\rho|y|^{2k} |\O|^2 ) +2(a+b-bk) \psi_\rho|y|^{2k} |\O|^2  +b(y\cdot \nabla )(\psi_\rho|y|^{2k} |\O|^2 ) +(V\cdot \nabla )(\psi_\rho|y|^{2k}|\O|^2 )}\hspace{.1in}\n \\
&&=2\hat{A}(y,s) \psi_\rho|y|^{2k} |\O|^2 +b |y|^{2k}|\O|^2 (y\cdot \nabla )\psi_\rho +|y|^{2k}|\O|^2 (V\cdot \nabla )\psi_\rho +2k |y|^{2k-1} |\O|^2 V_r \psi_\rho\n \\
&& \leq \frac12 (a+b-bk) \psi_\rho|y|^{2k}|\O|^2+ \frac{b}{2}  |\O|^2|y|^{2k}(y\cdot \nabla )\psi_\rho
+\frac12 (a+b-bk)\psi_\rho|y|^{2k} |\O|^2\n \\
 \eq
 for all $\rho \geq r_0$ and $s\in [0, S_0]$. Since $\frac{ b}{2}|y|^{2k} |\O|^2(y\cdot \nabla )\psi_\rho \leq 0$, we have the following differential inequality from (\ref{215}),
 \bb\label{216}
 \frac{\partial}{\partial s} f(y,s)+(a+b-bk) f(y,s)+b y\cdot \nabla f(y,s)+V\cdot \nabla f(y,s)\leq 0,
 \ee
 where we set
 $$
 f(y,s):= \psi_{r_0}|y|^{2k} |\O|^2.
 $$
 Let us define the space-time domain
 $$\mathcal{D}_{r_0}:=(\Bbb R^3\setminus  B(0, r_0)) \times (0, S_0).$$
 We have $\partial \mathcal{D}_{r_0}= \Gamma_1 \bigcup \Gamma_2 \bigcup \Gamma_3 \bigcup \Gamma_4$, where
\bqn &&\Gamma_1=\partial B(0, r_0)\times  (0, S_0), \quad \Gamma_2=\{ |y|=\infty\}\times  (0, S_0) \\
&& \Gamma_3= (\Bbb R^3\setminus  B(0, r_0))   \times \{s= S_0\} ,\quad \Gamma_4= (\Bbb R^3\setminus B(0,r_0))\times \{ s=0\}.
\eqn
The condition (\ref{decayiia}) implies that $f=0$ on $\Gamma_2$, and by construction of $\psi_{r_0}$ we have also
 $f(y,s)=0$ on $ \Gamma_2$.
Thus, the positive maximum of $f(y,s)$  on $\mathcal{\overline{D}}_{r_0} $, which exists due to assumption of the non-compactness of the support of $\O (y,s)$ in $\Bbb R^3\times [0, S_0]$, is attained only at $\Gamma_4$, and we have
 \bb\label{217}
 \sup_{|y|>r_0, s\in (0, s_0)}(\psi_{r_0}(y )|y|^{2k}|\O(y,s)|^2)=  \sup_{|y|>r_0}\psi_{r_0}(y)|y|^{2k} |\O_0(y)|^2 .
\ee
Moreover,   there exists $y_0 \in \Bbb R^3$ with $|y_0|>r_0$ such that the strict inequality
\bb\label{218}
\psi_{r_0}(y )|y|^{2k}|\O(y,s)|^2<  |\psi_{r_0}(y_0 )|y_0|^{2k}|\O_0(y_0)|^2= \sup_{|y|>r_0}\psi_{r_0}(y )|y|^{2k}|\O_0(y)|^2
\ee
holds for each $(y,s) \in \mathcal{D}_{r_0}\bigcap \Gamma_3 $.  Substituting $y=y_0, s=S_0$ into (\ref{218}), we have
\bb\label{219}
\psi_{r_0}(y_0 )|y_0|^{2k}|\O(y_0,S_0)|^2<   |\psi_{r_0}(y_0 )|y_0|^{2k}|\O_0(y_0)|^2.
\ee
Since $|\O(y_0,S_0)|= |\O_0(y_0)|>0$ by the periodicity in time, (\ref{219}) is absurd. \\
 \ \\
\noindent{\em \underline{(D) The case of (ii) with $a+b\geq 0, b>0 $:} }
Similarly to the proof (B) above  we introduce  $\overline{V} (y, s)= V(y, S_0-s)$,  $\bar{P} (y, s)= P(y, S_0-s)$ and $\Ob (y, s)= \O(y, S_0-s)$ for $0\leq s\leq S_0$ to derive (\ref{211}), and then we are reduced to the case of (C).
$\square$\\
 \ \\
\noindent{\bf Proof of Lemma 2.2 }
Similarly to (\ref{calculus}), from
 $$ V(y,s)=V(0,s)+\int_0 ^1  y\cdot V(\s y,s) d\s $$
 we obtain
 $|V(y,s)|\leq |V(0,s)|+ |y|\|\nabla V(s)\|_{L^\infty}\leq C(1+|y|)(\|\nabla
 V(s)\|_{L^\infty}+|V(0,s)|),$
 for each $s\in [0, S_0]$, and
  \bb\label{calb}
  \sup_{y\in \Bbb R^3} \frac{|V(y,s)|}{1+|y|} \leq C(\|\nabla
  V(s)\|_{L^\infty}+|V(0,s)|).
  \ee
  Given $R>0$, we introduce the radial cut-off function $\sigma_R(y)= 1-\psi_R (y)$, where $\psi_R (y)$ is the function $\psi_\rho (y)$ defined in (\ref{cut})
  with $\rho=\delta=R$.
 Let $\delta >0$. We  take $L^2(\Bbb R^3\times [0, S_0])$ inner product (\ref{vol}) by $\O (\delta +|\O|^2)^{\frac{q}{2}-1}
 \si$ to obtain
 \bq\label{newvorb}
&& (a+b)\int_0 ^{S_0} \int_{\Bbb R^3} |\O|^2(\delta +|\O|^2)^{\frac{q}{2}-1}\si dyds -\int_0 ^{S_0} \int_{\Bbb R^3} [(\O\cdot \nabla )V]\cdot\O (\delta +|\O|^2)^{\frac{q}{2}-1}\si dyds\n \\
&&\qquad = -\frac{1}{q}\int_0 ^{S_0} \int_{\Bbb R^3} \left[(( by+V)\cdot \nabla
)(\delta +|\O|^2)^{\frac{q}{2}}\right]\si dyds,\n \\
 \eq
 where we used the fact
  $$ \int_0 ^{S_0} \int_{\Bbb R^3} \frac{\partial \O}{\partial s} \cdot \O(\delta +|\O|^2)^{\frac{q}{2}-1}
 \si dyds = \frac{1}{q}\int_0 ^{S_0} \frac{\partial}{\partial s}  \int_{\Bbb R^3}(\delta +|\O|^2 )^{\frac{q}{2}}\si dy=0
 $$
 due to the periodicity.
 For fixed $\delta>0$ and $R>0$ the integrands in the right hand side of (\ref{newvorb}) are sufficiently smooth functions having the compact support, and one can integrate by part  them to obtain
 \bq\label{newvor1b}
&& (a+b)\int_0 ^{S_0}\int_{\Bbb R^3} |\O|^2(\delta +|\O|^2)^{\frac{q}{2}-1}\si dyds - \int_0 ^{S_0} \int_{\Bbb R^3}[ (\O\cdot \nabla )V]\cdot\O (\delta +|\O|^2)^{\frac{q}{2}-1}\si dyds\n \\
&&\qquad = \frac{3b}{q}\int_0 ^{S_0} \int_{\Bbb R^3}(\delta +|\O|^2)^{\frac{q}{2}}\si dx
+\frac{1}{q}\int_0 ^{S_0} \int_{\Bbb R^3} (\delta +|\O|^2)^{\frac{q}{2}} \left(( by+V)\cdot \nabla \right)\si dyds.\n \\
 \eq
 Passing $\delta \downarrow 0$ in (\ref{newvor1b}), using the dominated convergence theorem, we have
 \bq\label{selb}
 \lefteqn{\left(a+b- \frac{3b}{q}\right)\int_0 ^{S_0} \int_{\Bbb R^3} |\O|^q \si dyds-
\int_0 ^{S_0} \int_{\Bbb R^3} (\O \cdot \nabla) V\cdot \O  |\O|^{q-2}\si\, dyds}\n \\
&&=\frac{ b}{q}\int_0 ^{S_0}  \int_{\Bbb R^3} |\O|^q (y\cdot\nabla )\si \,
dx +\frac{1}{q} \int_0 ^{S_0} \int_{\Bbb R^3} |\O|^q (V\cdot\nabla )\si \, dyds\n \\
&&:=I+J.
 \eq
We estimate $I$ and $J$ easily as follows.
  $$
  |I|\leq  \frac{|b|}{ qR}\int_0 ^{S_0} \int_{\{R\leq |x|\leq 2R\}} |\O|^q
  |y||\nabla \s|\, dy
  \leq \frac{2|b|}{q}\|\nabla \s\|_{L^\infty}\|\O\|^q_{L^p(\{R\leq |y|\leq
  2R\}\times [0, S_0])}\to 0
  $$
   as $R\to \infty$.
   \bqn
    |J|&\leq&  \frac{1}{q R}\int_0 ^{S_0} \int_{\{R\leq |y|\leq 2R\}} |\O|^q
  |V||\nabla \s|\, dy
  \leq\frac{1+2R}{q R}\int_0 ^{S_0} \int_{\{R\leq |y|\leq 2R\}}
  \frac{|V(y,s)|}{1+|y|}|\O|^q|\nabla \s|\, dyds\n \\
  &\leq &\frac{C(1+2R)}{q R}\|\nabla \s\|_{L^\infty}\sup_{s\in [0, S_0]}(\|\nabla V(s)\|_{L^\infty}+|V(0,s)|) \|\O\|^q_{L^p(\{R\leq |y|\leq
  2R\}\times [0, S_0])}\to 0
  \eqn
  as $R\to \infty$, where we used (\ref{calb}). Therefore, passing
  $R\to \infty$ in (\ref{selb}), and using the dominated convergence theorem for the left hand side, we obtain,
  $$
  \left(a+b- \frac{3b}{q}\right)\int_0 ^{S_0}\int_{\Bbb R^3} |\O|^q  dyds=
\int_0 ^{S_0}\int_{\Bbb R^3} (\O \cdot \nabla) V\cdot \O  |\O|^{q-2}\, dyds,
  $$
  from which we deduce easily
  \bq\label{sel2b} -\sup_{s\in [0, S_0]}\|\nabla V\|_{L^\infty}\int_0 ^{S_0}\int_{\Bbb R^3} |\O|^q  dyds
&\leq& \left(a+b-\frac{3b}{q}\right)\int_0 ^{S_0}\int_{\Bbb R^3} |\O|^q  dyds\n \\
&\leq&\sup_{s\in [0, S_0]}\|\nabla V\|_{L^\infty}\int_0 ^{S_0}\int_{\Bbb R^3} |\O|^q  dyds.\n \\
  \eq
 Suppose there exists $(y_0, s_0)\in \Bbb R^3\times [0, S_0]$ such that  $\O(y_0, s_0)\neq 0$, then since $\O$ is a continuous function on $\Bbb R^3\times [0, S_0]$, one has $\int_0 ^{S_0}\int_{\Bbb R^3} |\O|^q  dyds> 0$, and we can divide (\ref{sel2b}) by
 $\int_0 ^{S_0}\int_{\Bbb R^3} |\O|^q  dyds$ to have
\bb\label{sel3b}
 -\sup_{s\in [0, S_0]}\|\nabla V(s)\|_{L^\infty} \leq
\left(a+b-\frac{3b}{q}\right)\leq \sup_{s\in [0, S_0]}\|\nabla V(s)\|_{L^\infty},
  \ee
 which holds for all $q\in (0, r)$ and for some $r>0$. Since $b\neq 0$, passing $q\downarrow 0$ in (\ref{sel3b}), we
 obtain desired contradiction. Therefore $\O(\cdot, s)=\mathrm{curl}\, V(\cdot, s)=0$ for all $s\in [0, S_0]$. This,
 together with $\mathrm{div}\,V(\cdot, s)=0$, provides us with the fact that $V(\cdot, s)=\nabla
 h(s)$ for a scalar harmonic function $h$ on $\Bbb R^3$.  If we impose the condition $|\nabla V(y,s)|=|\nabla \nabla h(y,s)|\to 0$
 as $|y|\to \infty$, then by the Liouville theorem for a harmonic function $\nabla V = 0$, and  $V(y,s)$ is independent of $y\in \Bbb R^3$, and
 periodic in $s$ with the period $S_0$.$\square$\\
 \ \\

\noindent{\bf Proof of Theorem 1.3 }
 We first note that
\bqn
  \int_{\Bbb R^3} |\xi| |\hat{v}(\xi,t)|d\xi &\leq& \left(\int_{\Bbb R^3} |\hat{v}(\xi,t)|^2 (1+|\xi|^2)^m d\xi\right)^{\frac12} \left(\int_{\Bbb R^3} \frac{|\xi|^2}{(1+|\xi|^2)^m} d\xi\right)^{\frac12}\n \\
  &\leq&  C\|v\|_{L^\infty([0, T); H^m (\Bbb R^3))} \left(\int_0 ^\infty \frac{r^4}{(1+r^{2m})} dr\right)^{\frac12}
  <\infty
  \eqn
  for $m>5/2$, and thus $\widehat{\nabla v}(\xi,t)\in L^1 (\Bbb R^3)$ for all $t\in [0, T)$. Hence, by the Riemann-Lebesgue lemma for the inverse Fourier transform we have
\bb\label{rl}|\nabla v(x,t)| \to 0 \quad \mbox{as $|x|\to \infty$} \quad \forall t\in [0, T).
\ee
Let us recall that the trajectory mapping $X_t(a)$ is defined by
\bb\label{230}
 \frac{\partial X_t(a)}{\partial s} = v(X_t(a),t), \quad X_0 (a)=a.
 \ee
Let  $\psi_\rho$ be the cut-off function defined in the proof of Theorem 1.2.  We define a moving cut-off by $\Psi_\rho(x,t)= \psi_\rho (X_t ^{-1} (x))$, which satisfies the transport equation,
 $$ \left\{\aligned  &\frac{\partial \Psi_\rho(x,t)}{\partial t}+(v\cdot \nabla) \Psi_\rho (x,t)=0,\\
 &\Psi_\rho (x,0)= \psi_\rho (x),\endaligned \right.
 $$
 and therefore
\bb\label{231}
\Psi_\rho (x,t)=\left\{ \aligned &1, \quad \mbox{if}\quad x\in \Bbb R^3\setminus X_t ( B(0, \rho+\delta)),\\
                          &0,     \quad   \mbox{if} \quad x\in X_t (B(0, \rho)),
                                \endaligned
                                \right.
                               \ee
 and
 $$0\leq \Psi_\rho (x,t)\leq 1 \quad\mbox{if}\quad x\in X_t ( B(0, \rho+\delta))\setminus X_t (B(0, \rho)).$$
Taking curl of (E), and multiplying it by $\o\Psi_\rho e^{-2\vare t}$, we have
 \bb\label{232}
  \frac{\partial }{\partial t}(e^{-2\vare t}|\o|^2\Psi_\rho  )+ 4\vare  e^{-2\vare t}| \o|^2 \Psi_\rho+(v\cdot \nabla ) (e^{-2\vare t}|\o|^2\Psi_\rho)=2\hat{\a}e^{-2\vare t}|\o |^2\Psi_\rho,
  \ee
 where $\hat{\a}(x,t)$ is defined by
 $$ \hat{\a}(x,t)=\left\{ \aligned & \frac{(\o(x,t))\cdot \nabla )v(x,t)\cdot \o(x,t)}{|\o(x,t)|^2}, \quad \mbox{if}\quad \o(x,t)\neq 0,\\
 & 0, \quad \mbox{if}\quad \o(x,t)= 0.
 \endaligned \right.
 $$
 We fix $t_0\in (0, T)$. By (\ref{rl}) there exists $R_1=R_1( v_0, t_0, \vare ) >0$ such that
 $$
 |\nabla v(x,t)|\leq 2\vare \quad \mbox{if $ |x|\geq R_1$ and $t\in [0, t_0]$}.
 $$
 For such $R_1$ we set $R_0=R_0 (v_0, t_0, \vare) $ as
 $$
R_0 = R_1+\int_0 ^{t_0} \|v(t)\|_{L^\infty} dt.
 $$
 Then, we have
 \bb\label{233}
 \min_{0\leq t\leq t_0} \mathrm{dist}(0, X_t (\partial B(0, R_0)))\geq R_1,
 \ee
 and therefore
\bb\label{234}
 |\nabla v(x,t)|\leq 2\vare\quad \mbox{if}\quad x\in \Bbb R^3\setminus X_t ( B(0, r_0))
\ee
 for all $t\in [0, t_0]$, and
\bb\label{235}
|\hat{\a}(x,t)| |\o(x,t)|^2\Psi_{R_0}(x,t) \leq |\nabla v(x,t)| |\o(x,t)|^2\Psi_{R_0}(x,t)\leq  2\vare |\o(x,t)|^2\Psi_{R_0}(x,t)
 \ee
 for all $(x,t)\in \Bbb R^3\times [0, t_0]$.
Substituting (\ref{235})  into (\ref{232}), and absorbing the vortex stretching term into the left hand side, one obtains
 \bb\label{236}
\frac{\partial}{\partial t}f(x,t)+ (v\cdot \nabla )f(x,t) \leq 0, \quad f(x,t):=e^{-2\vare t}|\o(x,t)|^2\Psi_{R_0}(x,t).
 \ee
 This implies that $f(X_t (a),t)\leq f(a,0)$, which is written in terms of original representation as
 \bb\label{237}
 e^{-2\vare t} |\o ( X_t (a), t)|^2 \psi_{R_0} (a)\leq |\o_0 (a)|^2 \psi_{R_0} (a),\quad \forall t\in [0, t_0]
 \ee
 where we used the fact $\Psi_{R_0} (X_t (a),t) =\psi_{R_0} ({X_t} ^{-1}(X_t (a)))=\psi_{R_0} (a)$.
 Since $\psi_{R_0}(a)>0$ for $|a|>R_0$, from (\ref{237}) we have (\ref{est}). $\square$

\section{ Remarks on the MHD system}
\setcounter{equation}{0}
The proof of the previous section applies well to the more general system such as the inviscid magnetohydrodynamic equations in $\Bbb R^3$,
\[
\mathrm{ (MHD) }
 \left\{ \aligned
 &\frac{\partial v}{\partial t} +(v\cdot \nabla )v =(b\cdot\nabla)b-\nabla (p +\frac12 |b|^2), \\
 &\frac{\partial b}{\partial t} +(v\cdot \nabla )b =(b \cdot \nabla
 )v,\\
 &\quad \textrm{div }\, v =\textrm{div }\, b= 0 ,\\
  &v(x,0)=v_0 (x), \quad b(x,0)=b_0 (x),
  \endaligned
  \right.
  \]
where $v=(v_1, v_2 , v_3 )$, $v_j =v_j (x, t)$, $j=1,\cdots,3$,
is the velocity of the flow, $p=p(x,t)$ is the scalar pressure,
$b=(b_1, b_2 , b_3 )$, $b_j =b_j (x, t)$, is the magnetic field,
and $v_0$, $b_0$ are the given initial velocity and magnetic field,
 satisfying div $v_0 =\mathrm{div}\, b_0= 0$, respectively. If we set $b=0$ in (MHD), then it reduces to (E).
  The scaling property of (MHD) is the same as (E).
 Given $\alpha \neq -1$, we  make the  self-similar
transform of (MHD), which is defined by the map $(v,b, p) \mapsto (V,B,P)$
given by
\begin{equation}
\label{31}
 v(x, t)=\frac {1}{(T-t)^{\frac{\a}{\a+1}}} V(y,s), \,
 b(x, t)=\frac {1}{(T-t)^{\frac{\a}{\a+1}}} B(y,s),\,
 p(x,t)=\frac {1}{(T-t)^{\frac{2\a}{\a+1}}} P(y,s),
\end{equation}
where
\begin{equation}
\label{32} y = \frac{x-x_*}{(T-t)^{\frac{1}{\a+1}}}, \quad s =  \log \left(\frac{T}{T-t}\right).
\end{equation}
Substituting (\ref{31})--(\ref{32}) into (MHD), we have the following system in terms of  $(V,B,P)$:
\begin{equation}\label{mhd} \left\{
\aligned
&\frac {\partial V}{\partial s}+ \frac{\a}{\a+1} V +\frac{1}{\a+1}(y \cdot \nabla)V + (V\cdot \nabla )V =
(B\cdot \nabla ) B-\nabla
(P+\frac12 |B|^2),\\
&\frac {\partial B}{\partial s}+ \frac{\a}{\a+1}B +\frac{1}{\a+1}(y \cdot \nabla)B + (V\cdot \nabla )B =(B\cdot \nabla ) V,\\
& \mathrm{div}\,V=0,\,\,\,\mathrm{div}\,B=0, \\
& V(y,0)=V_0 (y)=T^{\frac{\a}{\a+1}}v_0 (T^{\frac{1}{\a+1}}y),\, B(y,0)=B_0 (y)=T^{\frac{\a}{\a+1}}b_0 (T^{\frac{1}{\a+1}}y).
\endaligned \right.
\end{equation}
The following theorem shows the triviality of the time periodic solution of (\ref{mhd}), which implies that there exists no nontrivial discretely self-similar blowing up solutions to the system (MHD), if the blow-up profile $(V, B)$ satisfies suitable decay conditions near infinity as described in the theorem. This is a substantial improvement of the main result in \cite{cha3}.
\begin{thm}
Let $(V,B)\in C^1_sC^2_y (\Bbb R^{3+1})$ is a time periodic solution to (\ref{mhd}) with the period $S_0 \neq 0$,
  satisfying
  \bb
 |\nabla V( y,s)| +|B(y,s)|=o(1), \quad \forall s\in [0, S_0].
  \ee
  We assume furthermore that the vorticity $\O(y,s)$  and $B(y,s)$ decay depending $\a$ as follows.
  \begin{itemize}
\item[(i)] In the case  $\a <-1$: No extra condition.
\item[(ii)] In the case $-1< \a <0$:  There exists $k> \a+1$ such that
\bb
\ |\O(y,s)|=O(|y|^{-k})\quad \mbox{as $|y|\to \infty$}\quad \forall s\in [0, S_0].
\ee
\item[(iii)] In the case $\a \geq 0$:  There exist $k_1>\a+1$ and $k_2 >\a$ such that
\bb
\ |\O(y,s)|=O (|y|^{-k_1}), \quad |B(y,s)|=O(|y|^{-k_2}) \quad \mbox{as $|y|\to \infty$}\quad \forall s\in [0, S_0].
\ee
\end{itemize}
Then,  $V(y,s)=C(s), B(y,s)=0$ for all $(y,s)\in \Bbb R^3\times[0, S_0]$, where $C:[0, S_0] \to \Bbb R^3$ is a closed curve, $C(s)=C(s+S_0)$ for all $s\in [0, S_0)$.
\end{thm}
\noindent{\bf Proof } We first apply Lemma 2.1 and Lemma 2.2 to the second equations of (\ref{mhd}). Note that in $(V, B)$ of (\ref{mhd}) corresponds to
$(V, \O)$ of (\ref{vol}) respectively, and the coefficients $a+b$ and $b$ in (\ref{vol})  have roles of $\frac{\a}{\a+1}$ and
$\frac{1}{\a+1}$ respectively.  The following the proof of Lemma 2.1, we conclude that $B$ has compact support, and then by Lemma 2.2, we conclude that $B=0$. Now, since $B=0$ in (\ref{mhd}), the problem reduces to that of the Euler system (E).
We now apply  Theorem 1.1 to conclude that $V(y,s)=C(s)$, a closed curve in $\Bbb R^3$ with $C(s)=C(s+S_0)$. $\square$\\
\ \\
 $$ \mbox{\bf Acknowledgements } $$
 The author would like to thank to the anonymous referee for careful reading and constructive criticism.
 This research is supported partially by NRF
  Grants no.
 2006-0093854 and  no. 2009-0083521.

\end{document}